\documentclass[runningheads]{llncs}
\usepackage[utf8]{inputenc}

\usepackage{cite}
\usepackage{amsmath}
\usepackage{amsfonts}
\usepackage{amssymb}
\usepackage{authblk}
\usepackage{graphicx}
\usepackage{subcaption}
\usepackage[overload]{empheq}
\usepackage{mdwlist}
\usepackage[ruled,vlined]{algorithm2e}

\usepackage{pgfplots}
\usepackage{tikz}
\usepackage{float}
\usepackage{cleveref}
\usepackage[normalem]{ulem}

\usetikzlibrary{calc} 
\usetikzlibrary{shapes,arrows,positioning,shadows,trees}
\usetikzlibrary{decorations.pathreplacing}
\usetikzlibrary{positioning}
\usetikzlibrary{shapes.misc}
\usetikzlibrary{patterns}
\usetikzlibrary{3d}
\usetikzlibrary{spy}
\usetikzlibrary[patterns]
\usetikzlibrary{arrows,positioning,shapes} 
\usetikzlibrary{datavisualization}

\usepackage{tikzscale}
\usepgfplotslibrary{external} 
\begin{document}

\title{On quadrature rules for solving Partial Differential Equations using Neural Networks}
%
%
\author{Jon A. Rivera\inst{1,2} \and
Jamie M. Taylor\inst{2} \and
\'Angel J. Omella \inst{1} \and
David Pardo\inst{1,2,3}}
\authorrunning{J. A. Rivera et al.}
%
\institute{University of the Basque Country (UPV/EHU), Leioa, Spain \and
Basque Center for Applied Mathematics (BCAM), Bilbao, Spain \and Ikerbasque (Basque Foundation for Sciences), Bilbao, Spain}

\maketitle              

\begin{abstract}
Neural Networks have been widely used to solve Partial Differential Equations. These methods require to approximate definite integrals using quadrature rules. Here, we illustrate via 1D numerical examples the quadrature problems that may arise in these applications and propose different alternatives to overcome them, namely: Monte Carlo methods, adaptive integration, polynomial approximations of the Neural Network output, and the inclusion of regularization terms in the loss. We also discuss the advantages and limitations of each proposed alternative. We advocate the use of Monte Carlo methods for high dimensions (above 3 or 4), and adaptive integration or polynomial approximations for low dimensions (3 or below). The use of regularization terms is a mathematically elegant alternative that is valid for any spacial dimension, however, it requires certain regularity assumptions on the solution and complex mathematical analysis when dealing with sophisticated Neural Networks. 

\end{abstract}


\section{Introduction}

In the last years, the use of Deep Learning (DL) methods has grown exponentially in multiple areas, including self-driving cars \cite{DL_cars_1,DL_cars_2}, speech recognition \cite{DL_speech_1,DL_speech_2}, and healthcare \cite{DL_health_1,DL_health_2}. Similarly, the use of DL algorithms has become popular for solving Partial Differential Equations (PDEs) -- see, e.g., \cite{DL_PDE_1,DL_PDE_2,DL_PDE_3,DL_PDE_4, pinn_BC, DL_me,muga,Maciej}.

DL techniques present several advantages and limitations with respect to traditional PDE solvers based on Finite Elements (FE) \cite{book_fem},  Finite Differences (FD) \cite{book_FD}, or Isogeometric Analysis (IGA) \cite{IGA_1}. Among the advantages of DL, we encounter the nonnecessity of generating a grid. In general, DL uses a dataset in which each datum is independent from others. In contrast, in the linear system that produces the FEM method there exists a connectivity between the nodes of the mesh. The independence of the data in DL allows the parallelization into GPUs for fast computations. Furthermore, DL provides the possibility of solving certain problems that cannot be solved via traditional methods, like high-dimensional PDEs \cite{DL_HD_1,DL_HD_2}, some fractional PDEs \cite{DL_Fr_1,DL_Fr_2}, and multiple nonlinear PDEs \cite{DL_NL_1,DL_NL_2}.

However, DL also presents limitations when solving PDEs. For example, in \cite{wang2020pinns} they show that Fully-Connected Neural Networks suffer from spectral bias \cite{spectral_bias} and exhibit different convergence rates for each loss component. In addition, the convergence of the method is often assumed (see, e.g., \cite{mishra2020estimates}) since it cannot be rigorously guaranteed due to the non-convexity of the loss function. Another notorious problem is due to quadrature errors. In traditional mesh-based methods, such as FEM, we first select an approximating solution space, and then we compute the necessary integrals over each element required to produce the stiffness matrix. With DL, we first set a quadrature rule and then construct the approximated function. As we will show throughout this work, this process may lead to large quadrature errors due to the unknown form of the integrand. This can also be interpreted as a form of overfitting \cite{over_1,over_2} that applies to the PDE constraint (e.g., $u"=0$) rather than to the PDE solution (e.g., $u$). This may have disastrous consequences since it may result in approximate solutions that are far away from the exact ones at all points.

There exist different methods to approximate definite integrals by discrete data. The most common ones used in DL are based on Monte Carlo integration. These methods compute definite integrals using randomly sampled points from a given distribution. Monte Carlo integration is suitable for high-dimensional integrals. However, for low-dimensional integrals (1D, 2D, and 3D), convergence is slow in terms of the number of integration points in comparison to other quadrature rules. This produces elevated computational costs. Examples of existing works that follow a Monte Carlo approach are \cite{DL_Ritz} using a Deep Ritz Method (DRM)  and \cite{DL_GM} using a Deep Galerkin Method (DGM). From the practical point of view, at each iteration of the Stochastic Gradient Descent (SGD), they consider a mini-batch of randomly selected points to discretize the domain.

Another existing method to compute integrals in DNNs is the so-called \textit{automatic integration} \cite{auto_int}. In this method, the author approximates the integrand by its high order Taylor series expansion around a given point within the integration domain. Then, the integrals are computed analytically. The derivatives needed in the Taylor series expansion are computed via automatic differentiation (a.k.a \textit{autodiff}) \cite{autodiff_2}. Since we are using the information of the derivatives locally, we are doomed to perform an overfitting on those conditions.

Another alternative is to use adaptive integration methods. In the Deep Least Square (DLS) method \cite{DL_LS}, they use an adaptive mid-point quadrature rule using local error indicators. These indicators are based on the value of the residual at randomly selected points, and the resulting quadrature error is unclear. In Variational Neural Networks (VarNets) \cite{DL_VarNet}, they use a Gauss quadrature rule to evaluate their integrals. As error indicators, they use the strong form of the residual evaluated over a set of random points that follow a uniform distribution.

One can also use Gauss-type quadrature rules to evaluate integrals, as they do in Variational Physics-Informed Neural Networks (VPINNs) \cite{DL_PINN_2}. One limitation of this method is that we cannot select an adequate quadrature order because we are unaware of the solution properties. In addition, as we are using fixed quadrature points, the chances of performing overfitting are high.

This works analyzes the problems associated with quadrature rules in DL methods when solving PDEs. As they remark in \cite{DL_PINN_2}, there is no proper quadrature rule in the literature developed for integrals of Deep Neural Networks (DNNs). In particular, we illustrate the disastrous consequences that may appear as a result of an inadequate selection of a quadrature rule. For that, we use a simple one-dimensional (1D) Laplace problem with the Deep Ritz Method (DRM) \cite{DL_Ritz}. We also show a theoretical example of the Deep Least Squares method. Then, we propose several alternatives to overcome the quadrature problems: Monte Carlo integration, exact (Gaussian) integration of piecewise-polynomial approximations, adaptive integration, and the use of regularizers. We discuss their advantages and limitations and develop adequate regularizers for certain problems using similar ideas to those presented in \cite{mishra2020estimates}. Moreover, we illustrate with a numerical example the different approximated solutions obtained with the different integration strategies.

The remainder of this article is as follows. Section \ref{sec:model_problem} presents our selected model problem. Section  \ref{sec:loss_func} defines the loss functions used in this work for solving PDEs. Section \ref{sec:NN} describes our Neural Network (NN) and explains some critical implementation aspects. Section \ref{sec:qua_rules} illustrates the quadrature problems via two numerical examples. Section \ref{sec:integral_approx} proposes several methods to overcome the quadrature problems and Section \ref{sec:num_res} shows numerical results corresponding to the proposed integration strategies. Finally, Section \ref{sec:conclu} summarizes the main findings and possible future lines of research in the area.


\section{Model problem} \label{sec:model_problem}

Let $\Omega \subset R^d$ be a computational domain, and $\Gamma_D$ and $\Gamma_N$ two disjoint sections of its boundary, where $\Gamma_D \cup \Gamma_N = \partial \Omega $ and the subscripts $D$ and $N$ denote the Dirichlet and Neumann bounds, respectively. We consider the following boundary value problem:

\begin{subequations}
\begin{align}[left=\empheqlbrace]
     -\nabla \cdot (\sigma \nabla u) & = f \qquad x \in \Omega, \label{eq:original_homog_problem_poisson_line1} \\
     u & = 0 \qquad x \in \Gamma_D, \label{eq:original_homog_problem_poisson_line2} \\
     (\sigma \nabla u) \cdot \bf{n}  & = g \qquad x \in \Gamma_N. \label{eq:original_homog_problem_poisson_line3} 
\end{align}
\label{eq:original_homog_problem_poisson} 
\end{subequations}
In the above, $\bf{n}$ is the unit normal outward (to the domain) vector, and we assume $\sigma(x) > 0$ for every $x$ and the usual regularity assumptions, namely, $\sigma \in L^\infty (\Omega)$, $f \in L^2 (\Omega)$, $g \in H^{1/2} (\Gamma_N)$, and $u \in V= H^{1}_{0}(\Omega) =\{v \in H^{1}(\Omega) \textup{  and  } v|_{\Gamma_D}= 0 \}$, where $H^{1}(\Omega) = \{ v \in L^{2}(\Omega), \nabla v \in L^{2}(\Omega)  \}$.

In this work we solve two different model problems to illustrate our numerical results. The solution to model problem 1 is $u(x)=x^{0.7}$ and it satisfies Eq. \eqref{eq:model_problem_1}.
\begin{equation}
\left\{
\begin{array}{rrcl}
     -u''(x) & = & 0.21x^{-1.3} & \qquad x \in (0,10),
 	 \\ u (0) & = & 0,
  	 \\ u'(10)  & = & \frac{0.7}{10^{0.3}}.
\end{array}
\right.
\label{eq:model_problem_1}
\end{equation}

The solution to model problem 2 is $u(x)=x^{2}$ and it satisfies Eq. \eqref{eq:model_problem_2}.
\begin{equation}
\left\{
\begin{array}{rrcl}
     u''(x) & = & 2 & \qquad x \in (0,10),
  \\ u (0) & = & 0,
  \\ u'(10)  & = & 20.
\end{array}
\right.
\label{eq:model_problem_2}
\end{equation}


\section{Loss functions} \label{sec:loss_func}
We introduce the following standard $L^2(\Omega)$ inner products:
\begin{equation}
(u,v) = \int_{\Omega} u \: v \qquad \textup{and} \qquad (g,v)_{\Gamma_N} = \int_{\Gamma_N} g \: v.
\label{eq:notation_1}
\end{equation} 

In the following, we consider two methods: the Ritz Method \cite{Ritz1909}, and the Least Squares Method \cite{LSM_theory}.

\subsection{Ritz Method}\label{sec:ritz}

Multiplying the PDE from Eq. \eqref{eq:original_homog_problem_poisson_line1} by a test function $v \in V$ (where $V=H^1_0(\Omega)$), integrating by parts and incorporating the boundary conditions, we arrive at the variational formulation:
\begin{equation} 
 \textup{Find  } u \in V \textup{  such that  } (\sigma \nabla u,\nabla v)=(f,v) + (g,v)_{\Gamma_N} \qquad \forall v \in V \label{eq:Ritz_variational}.
\end{equation}
It is easy to prove that the original (strong) and variational (weak) formulations are equivalent (see, for instance, \cite{Claes}).

To introduce the Ritz method, we define the energy function $\mathcal{F}_{R}: V \longrightarrow \mathbb{R}$ given by
\begin{equation}
\mathcal{F}_{R}(v) = \frac{1}{2} (\sigma \nabla v,\nabla v) - (f,v) - (g,v)_{\Gamma_N}.
\label{eq:Ritz_funct}
\end{equation}

As proved in \cite{Claes}, problem  \eqref{eq:Ritz_variational} is equivalent to the following energy minimization problem:
\begin{equation} 
 u = \underset{v \in V}{\textup{ arg min }} \mathcal{F}_{R}(v), \label{eq:Ritz_min}
\end{equation}

\subsection{Least Squares (LS) Method}

Reordering the terms of Eq. \eqref{eq:original_homog_problem_poisson_line1} and \eqref{eq:original_homog_problem_poisson_line3}, we define:
\begin{subequations}
\begin{align}
{\cal G}u := \nabla \cdot (\sigma \nabla u) + f \qquad x \in \Omega, \label{eq:LSM_1} \\
{\cal B}u := (\sigma \nabla u) \cdot \textbf{n} - g \qquad x \in \Gamma_N. \label{eq:LSM_2}
\end{align}
\label{eq:LSM} 
\end{subequations}
To introduce the Least Squares method, we define the function $\mathcal{F}_{LS}: V \longrightarrow \mathbb{R}$
\begin{equation}
\mathcal{F}_{LS}(v) = |({\cal G}v, {\cal G}v)| +  |({\cal B}v, {\cal B}v)_{\Gamma_N}|.
\label{eq:LSM_funct}
\end{equation}

We want to minimize the function $\mathcal{F}_{LS}(v)$ subject to the essential (Dirichlet) BCs. We often find the minimum by taking the derivative equal to zero and ending up with a linear system of equations. In the context of DL, we can simply introduce the above loss function $\mathcal{F}_{LS}(v)$ directly in our NN. Therefore, we want to find 

\begin{equation} 
 u = \underset{v \in V}{\textup{ arg min }} \mathcal{F}_{LS}(v). \label{eq:LSM_min}
\end{equation}


\section{Neural Network Implementation} \label{sec:NN}

We train a Neural Network, named $u_{NN}(x;\theta)$, with the following architecture. In this first part, we define the trainable part of our NN, with learnable parameters $\theta$. We call it $u_{\theta}$. It is composed by:

\begin{enumerate}
\item An input layer. This layer receives the data in the form of a $n \times d $ vector, where $n$ is the number of samples and $d$ is the dimension of the data of the problem.
\item $l$ hidden dense layers with $n_i$ neurons at each layer $i \in \{ 1, \cdots, l \}$ and a \textit{sigmoid} activation function. In our particular examples, we select $l=1$.
\item An output layer that delivers $u_{\theta}$.
\suspend{enumerate}

Now, we add more layers to our scheme in order to implement Equations \eqref{eq:Ritz_funct} or \eqref{eq:LSM_funct}. For that, we introduce:

\resume{enumerate}
\item A non-trainable layer to impose the Dirichlet boundary conditions. For that, we select a function $\phi(x)$ that satisfies the Dirichlet conditions of the problem and its value is nonzero everywhere else \cite{BC_implement}. In this work, we select the following $\phi(x)$ functions for 1D problems in the interval $\Omega=[a,b]$:

\begin{equation}
\phi(x) = \prod_{x_D \in \Gamma_D} (x-x_D).
\end{equation}

Then, we generate a new output of the NN: $u_{NN}(x;\theta)=\phi(x)u_{\theta}(x)$ that strongly imposes the homogeneous Dirichlet boundary conditions.

\item A non-trainable layer to compute the loss function $\mathcal{F}_{R}$ or $\mathcal{F}_{LS}$ following Eqs. \eqref{eq:Ritz_funct} or \eqref{eq:LSM_funct}. Within this layer, we evaluate the integrals and the derivatives. We consider different quadrature rules, being the quadrature points part of the input data of our NN, along with the physical points of the domain. For computing the derivatives, we use automatic differentiation, except in some specific cases, where we employ finite differences. These cases are explicitly indicated throughout the text.
\end{enumerate}

\begin{figure}[!htp]
\centering
\includegraphics[width=1\textwidth]{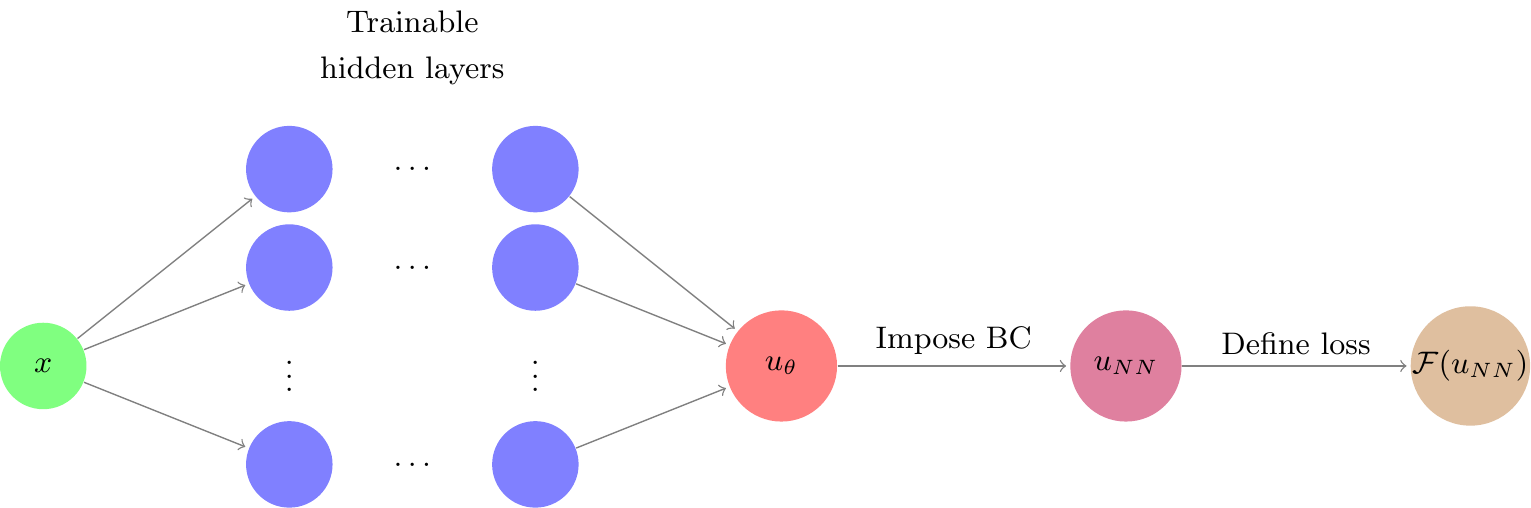}
	\caption{Our implementation scheme.}
	\label{fig:model_scheme}
\end{figure}

Figure \ref{fig:model_scheme} shows a schematic graph of the described NN architecture. Our software is developed in Python and we use the library \textit{Tensorflow 2.0}.

To train the NN, we replace in Equations \eqref{eq:Ritz_min} or \eqref{eq:LSM_min} the search space $V$ by the learnable parameters $\theta$ included in our NN. The result of the minimization is a function $u_{NN}(x, \tilde{\theta})$, where $\tilde{\theta}$ are the optimal learnable parameters encountered as a result of the training. For simplicity, in the following we abuse notation and use the symbol $u_{NN}$ to denote also the solution $u_{NN}(x, \tilde{\theta})$ of our minimization problem.

\section{Quadrature rules}\label{sec:qua_rules}

We approximate our integrals from Eqs. \eqref{eq:Ritz_funct} and \eqref{eq:LSM_funct} using a quadrature rule of the form

\begin{equation}
\int_a^b f(x) dx \approx \sum_{i=0}^{n} \omega_i f(x_i),
\label{eq:quadra_rule}
\end{equation}
where $\omega_i$ are the weights and $x_i$ are the quadrature points. Examples of quadrature rules that follow the above formula include trapezoidal rule and Gaussian quadrature rules \cite{book_Parviz}. We classify these quadrature rules into two groups: (1) those that only employ points from the interior of the interval; and (2) those that evaluate the solution at one extreme point or more (a or b). Integration rules within the later group (e.g., the trapezoidal rule) are inadequate for our minimization problems because the integrand can be infinite at the boundary points in the case of singular solutions. Thus, we focus on quadrature rules that only evaluate the solution at interior points of the domain, with a special focus on Gaussian quadrature rules.

\subsection{Illustration of quadrature problems in Neural Networks}

\subsubsection{Ritz method.}

We consider the two model problems from Section \ref{sec:model_problem}. We approximate $u(x)$ using the Ritz method. Thus, we search for a NN that minimizes the loss functional given by Eq. \eqref{eq:Ritz_funct}. Our NN has one hidden layer with $10$ neurons ($31$ trainable weights). We use automatic differentiation to compute the derivatives and a three-point Gaussian quadrature rule to approximate the integrals within each element. We select the Stochastic Gradient Descent (SGD) optimizer. For model problem 1, we discretize our domain with four equal-size elements and execute $40,000$ iterations during the optimization process. For model problem 2, we discretize our domain with ten equal-size elements and execute $200,000$ iterations during the optimization process.

\begin{figure}
\centering
\begin{subfigure}{.5\textwidth}
  \centering
  \includegraphics[width=0.9\linewidth]{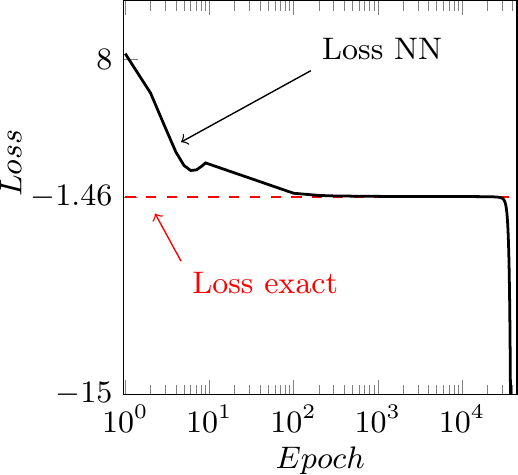}
  \caption{Model problem 1.}
  \label{fig:Ritz_model_problem_1_loss}
\end{subfigure}%
\begin{subfigure}{.5\textwidth}
  \centering
  \includegraphics[width=0.9\linewidth]{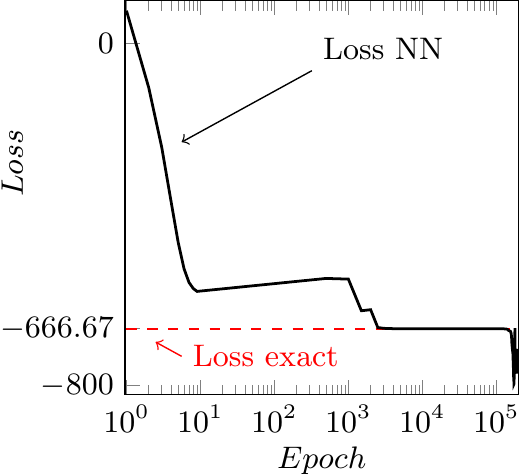}
  \caption{Model problem 2.}
  \label{fig:Ritz_model_problem_2_loss}
\end{subfigure}
\caption{Loss evolution of the training process for our two model problems.}
\label{fig:Ritz_model_problem_loss}
\end{figure}

Figures \ref{fig:Ritz_model_problem_1_loss} and \ref{fig:Ritz_model_problem_2_loss} describe the loss evolution of the training process. We obtain a lower loss than the best possible one (corresponding to the exact solution), which is due to high quadrature errors.


\begin{figure}[!htp]
\centering
	\subcaptionbox{Exact and approximate solutions.\label{fig:Ritz_model_problem_1a}}{%
	 \includegraphics[width=.42\linewidth]{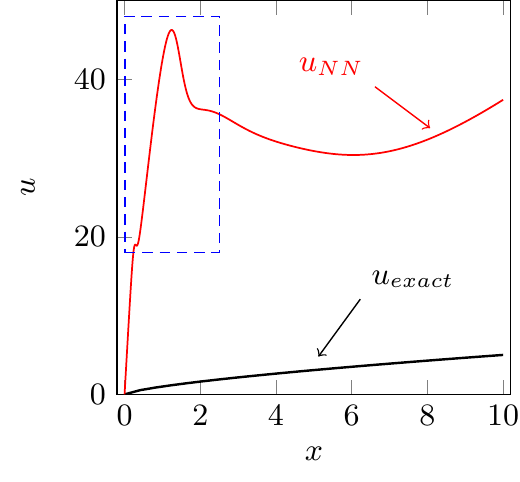}} 
	\hspace{1cm}
	\subcaptionbox{Approximate solution in the interval $[0,2.5]$. The Gauss quadrature points corresponding to the first element are indicated in blue.\label{fig:Ritz_model_problem_1b}}{%
	 \includegraphics[width=.4\linewidth]{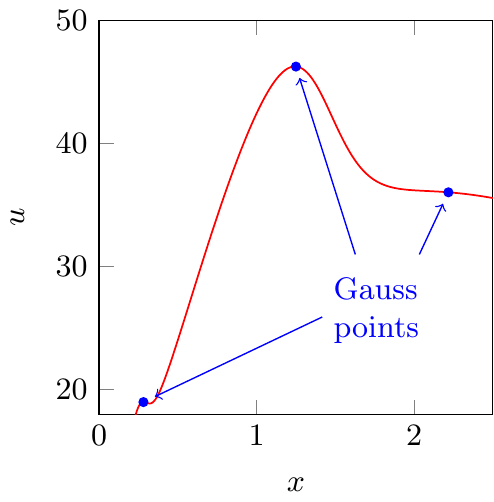}}
	\caption{Exact vs approximate Ritz method solutions of model problem 1 using four elements for evaluation of $\mathcal{F}_{R}(v)$ and a NN with 31 weights.}
	\label{fig:Ritz_model_problem_1}
\end{figure}

\begin{figure}[!htp]
\centering
	\subcaptionbox{Exact and approximate solutions.\label{fig:Ritz_model_problem_2a}}{%
	\includegraphics[width=.415\linewidth]{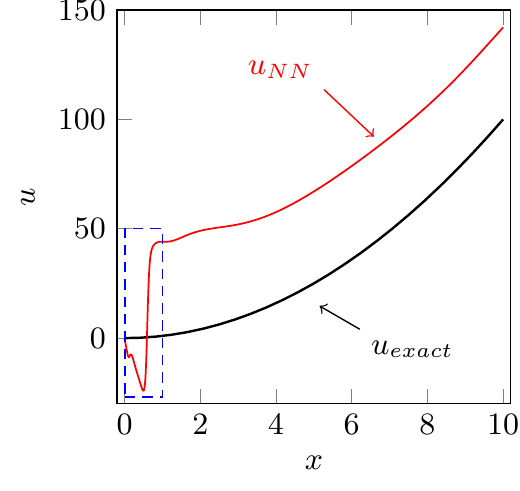}} 
	\hspace{1cm}
	\subcaptionbox{Approximate solution in the interval $[0,1]$. The Gauss quadrature points corresponding to the first element are indicated in blue.\label{fig:Ritz_model_problem_2b}}{%
	\includegraphics[width=.42\linewidth]{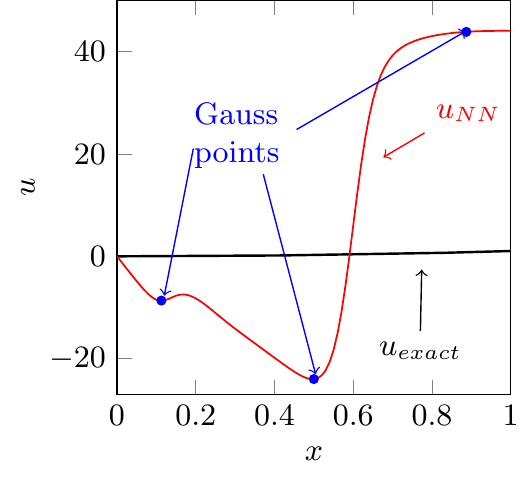}}
	\caption{Exact vs approximate Ritz method solutions of model problem 2 using ten elements for evaluation of $\mathcal{F}_{R}(v)$ and a NN with 31 weights.}
	\label{fig:Ritz_model_problem_2}
\end{figure}

Figures \ref{fig:Ritz_model_problem_1} and \ref{fig:Ritz_model_problem_2} compare the approximate and exact solutions. We observe a disastrous NN approximation due to quadrature errors. Figures \ref{fig:Ritz_model_problem_1b} and \ref{fig:Ritz_model_problem_2b} show that the gradient is (almost) zero at the training (quadrature) points of the first interval. This value minimizes the numerical aproximation of the integral of $(\nabla u)^2$. At the same time, the approximated solution is large in the first interval, thus maximizing the term $(f,u)$, and minimizing the total loss:

\begin{equation}
\notag
\mathcal{F}_{R}(u_{NN}) = \frac{1}{2} \underbrace{( \nabla u_{NN},\nabla u_{NN})}_{ \sum_{q_i} \omega_i (\nabla u_{NN})^2 \approx 0} - \underbrace{(f,u_{NN})}_{\approx \infty} - (g,u_{NN})_{\Gamma_N} \simeq - \infty
\notag
\end{equation}

The described quadrature errors can be interpreted as overfitting over a condition of the problem given by the energy of the solution.

\subsubsection{Least Squares Method.}

\begin{figure}[!htp]
\centering
  \includegraphics[width=.6\linewidth]{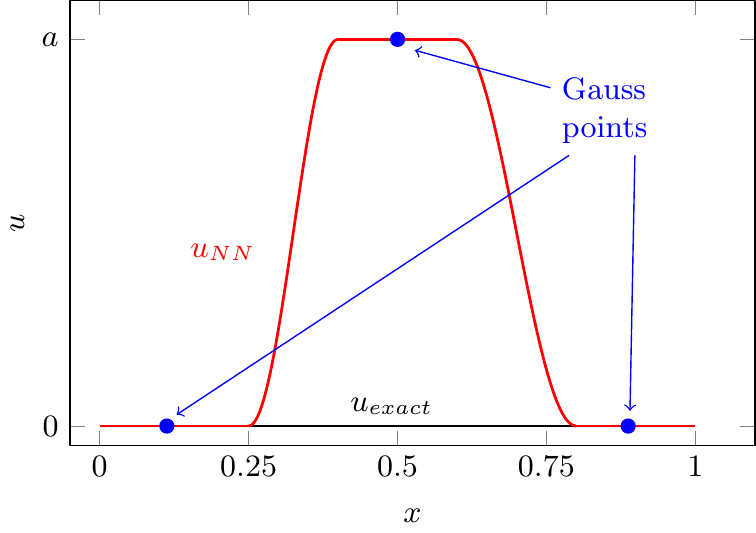}
	\caption{Exact ($u_{exact}=0$) and approximated solution of problem given by Eq. \eqref{eq:LSM_example} and solved with the LS method.}
	\label{fig:figure_LSM}
\end{figure}

We now consider the following one-dimensional problem:  

\begin{equation}
\left\{
\begin{array}{rrcl}
-u''(x)& = & 0 & \qquad x \in (0,1), \\ 
u(0) = u'(1) & = & 0, & \\ 
\end{array}
\right.
\label{eq:LSM_example}
\end{equation}
where the exact solution is $u(x) = 0$. We can easily construct an approximating function $u_{NN}$ that satisfies Eq. \eqref{eq:LSM_example} at the three considered Gaussian points and minimizes Eq. \eqref{eq:LSM_funct}, while still being a poor approximation of the exact solution due to quadrature errors. Figure \ref{fig:figure_LSM} shows an example.

\section{Integral approximation} \label{sec:integral_approx}

We now describe four different methods to improve the integral approximations.

\subsection{Monte Carlo integration}

We consider the following Monte Carlo integral approximation over a set of points $x_i \in (a,b)$, 

\begin{equation}
\int_{a}^{b} f(x) dx \approx \frac{(b-a)}{n} \sum_{i=1}^{n} f(x_{i}), \hspace{0.5cm} x_i \in (a,b) \hspace{0.2cm} \forall i=\{1,  \cdots, n \}
\end{equation}
In the above, points $x_i$ are randomly selected \cite{MC_2}. While this method is useful for higher-dimensional integrals, for lower dimensions (1D, 2D, 3D) the computational cost is high since the value of the integral approximation converges as $1/\sqrt{n}$ \cite{MC_1}.

\subsection{Piecewise-polynomial approximation}

\begin{figure}[!htp]
\centering
  \includegraphics[width=.75\linewidth]{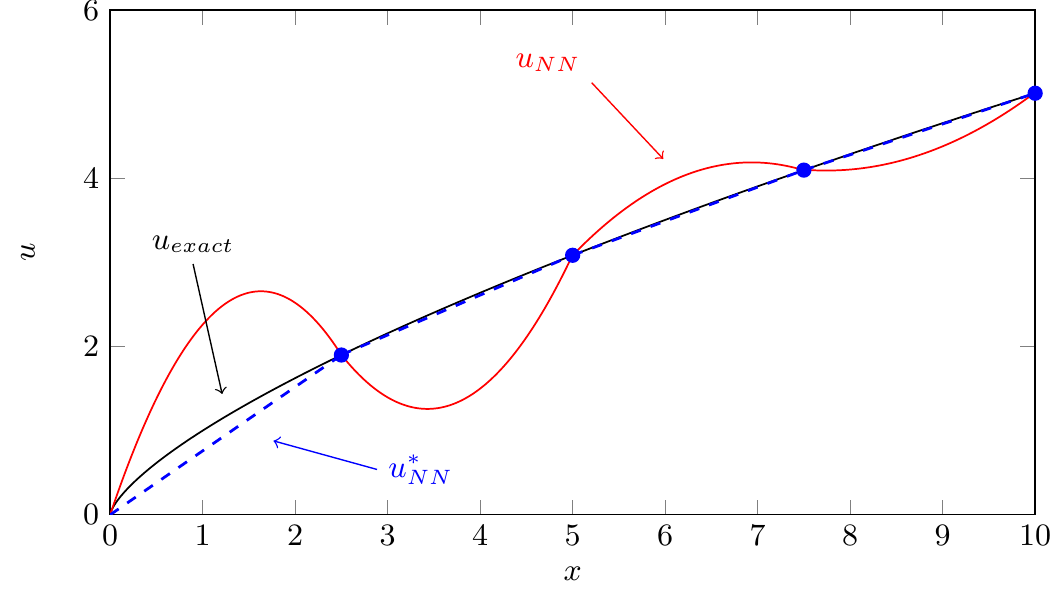}
	\caption{Neural Network approximation $u_{NN}$ and its piecewise-linear approximation $u^{*}_{NN}$.}
	\label{fig:linear_approx}
\end{figure}

We replace the original NN $u_{NN}$ by a piecewise-polynomial approximation $u^{*}_{NN}$, which can be exactly differentiated (e.q., via finite differences) and integrated (via a Gaussian quadrature rule). Figure \ref{fig:linear_approx} shows an example of the aforementioned case when we train a NN using four elements with linear approximations within each element.

This method controls quadrature errors. However, it is inadequate for high-dimensional problems as we need a mesh that is difficult to implement and integration becomes time consuming.

\subsection{Adaptive integration} 

\begin{figure}[!htp]
\centering
	\subcaptionbox{Training set\label{fig:adaptive_final_training}}{%
	  \includegraphics[width=.7\linewidth]{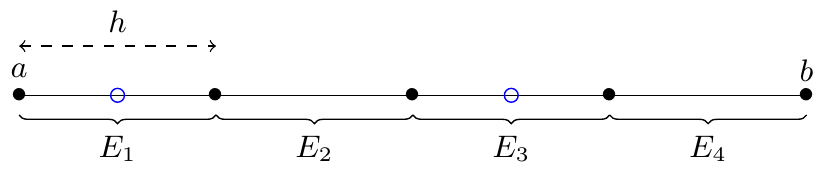}} 
	\subcaptionbox{Validation set\label{fig:adaptive_final_validation}}{%
	  \includegraphics[width=.7\linewidth]{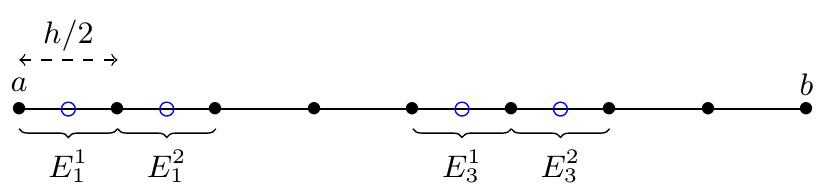}}
	\caption{Black points (dots) correspond to the original (a) training/(b) validation partitions and blue points (circles) are the points added by the refinement performed in the first and third elements.}
	\label{fig:adaptive_train_val}
\end{figure}

We first consider a training dataset over the interval $(a,b)$, and a validation dataset built from a global \textit{h}-refinement of the training dataset. Figure \ref{fig:adaptive_train_val} shows an example of a training and the corresponding validation datasets. Then, for each element of the training mesh (e.g., $E_1$ in Figure \ref{fig:adaptive_train_val}), we compare the numerical integral over that element vs the sum of the integrals over the two corresponding elements on the validation dataset (in our case, $E_1^1 + E_1^2$). If the integral values differ by more than a stipulated tolerance, we \textit{h}-refine the training element, and we upgrade the validation dataset so it is built as a global \textit{h}-refinement of the training dataset. This process is described in Algorithm \ref{alg:adaptive}. Figure \ref{fig:adaptive_train_val} also shows a training and the corresponding validation dataset after refining the first and third elements.

\begin{algorithm}[!htp]
\SetAlgoLined
 Generate a training dataset\;
 Generate the corresponding validation dataset\;
 Set tolerance $\epsilon$ and maximum iteration number $i_{max}$\;
 \While{$i < i_{max}$}{
 	\For{$j=1, \cdots, n \textup{ (number of elements)}$}{
  Compute integral values $I_j$ over the training dataset of elements $E_j$\;
  Compute integral values $I_j^1$,$I_j^2$ over the validation dataset of elements $E_j^1$,$E_j^2$ \;
  \eIf{$(I_j^1 + I_j^2) - I_j < \epsilon$}{
   $h$-refine the $E_j$-th element of the training set\;
   $h$-refine the $E_j^1$-th and $E_j^2$-th elements of the validation set\;
   }
   {
   continue\;
  }
  }
  $i = i + 1$\;
 }
 \caption{Adaptive integration method}
 \label{alg:adaptive}
\end{algorithm}

We are able to control the quadrature errors by adding new quadrature points to the training dataset. However, the simplest way to implement such method is by using meshes, which posses a limitation on high-dimensional integrals. As an alternative to generating a mesh, one can randomly add points to the training set. This entails difficulties when designing an adaptive algorithm.

In the same way that we propose an \textit{h}-adaptive method, we can also work with \textit{p}-adaptivity \cite{p_adap} or a combination of them (e.g., \textit{hp}-adaptivity \cite{hp_adap}).

\subsection{Regularization methods}\label{subsecHeuristic}

We now introduce a problem-specific regularizer designed to control the quadrature error.

In a one-dimensional setting, we consider the integral functional $\mathcal{F}_{R}$ as given in \eqref{eq:Ritz_funct}, and its approximation via a midpoint rule, $\hat{\mathcal{F}}_{R}$, given by 
\begin{equation}\label{eqIntegral}
\begin{split}
\hat{\mathcal{F}}_{R}(u)=&\frac{b-a}{N}\sum\limits_{j=1}^N \left(\frac{1}{2}|u'(x_j)|^2-f(x_j)u(x_j)\right)+g(b)u(b)+g(a)u(a). 
\end{split}
\end{equation}
We note that $g=0$, except where the Neumann condition is imposed. While we focus on the Ritz method, a similar heuristic can be applied to the Least Squares method. 

We introduce a function $\mathcal{R}$ that depends on the learnable parameters $\theta$ of a given neural network $u_{NN}$, such that for any neural network with a given architecture,
\begin{equation}\label{eqHeurErr}
|\mathcal{F}_{R}(u_{NN})-\hat{\mathcal{F}}_{R}(u_{NN})|<\mathcal{R}(\theta).
\end{equation}
If we then consider a loss function $\mathcal{L}$ given by 
\begin{equation}\label{eqRegLoss}
\mathcal{L}(\theta)=\hat{\mathcal{F}}_{R}(u_{NN})+\mathcal{R}(\theta),
\end{equation}
then we may be able to improve the approximation of the quadrature rule, as the loss contains a term that by design controls the quadrature error.

For simplicity, we consider only the case of a single-layer network with a one-dimensional input, and the midpoint rule for calculating the integral over a uniform partition of $(a,b)$. We consider a mid-point rule as in \eqref{eqIntegral}, and define the interval length $\delta=\frac{b-a}{N}$ and intervals $I_j=\left(\frac{\delta}{2}+x_j,x_j+\frac{\delta}{2}\right)$. We estimate the error of the midpoint rule to integrate $F$ as 

\begin{equation}\label{eqMidpError}
\begin{split}
\left|\int_a^b F(x)\,dx - \sum\limits_{j=1}^NF(x_j)\delta\right|=& \left|\sum\limits_{i=1}^N\int_{x_j-\frac{\delta}{2}}^{x_j+\frac{\delta}{2}}(F(x)-F(x_j))\,dx\right|\\
\leq & \sum\limits_{i=1}^N\int_{x_j-\frac{\delta}{2}}^{x_j+\frac{\delta}{2}}\left| F(x)-F(x_j)\right|\,dx\\
\leq & \sum\limits_{i=1}^N\max\limits_{t\in I_j}|F'(t)|\int_{x_j-\frac{\delta}{2}}^{x_j+\frac{\delta}{2}} |x-x_j|\,dx\\
=& \frac{\delta^2}{4}\sum\limits_{i=1}^N\max\limits_{t\in I_j}|F'(t)|.
\end{split}
\end{equation}

This estimate scales as $\mathcal{O}(\frac{1}{N})$ for fixed $F$, and thus, for a large number of integration points, we expect the estimate to be sufficiently accurate and to avoid ``overdamping" of the loss. 

With \eqref{eqMidpError} in mind, we estimate the local Lipschitz constants of the integrand as in \eqref{eq:Ritz_funct}. The numerical estimation of the Lipschitz constants of NNs has attracted attention, as they form a way of of estimating the generalizability of a neural network, and have been used in the training process as a way to encourage accurate generalization \cite{fazlyab2019efficient,gouk2021regularisation,scaman2018lipschitz}. As we are dealing with loss functions that involve {\it derivatives} of the neural network, we however need estimates of higher order derivatives of $u_{NN}$. The approach that we employ is similar in spirit to the work of \cite{mishra2020estimates} for obtaining {\it a posteriori} error estimates in PINNs. 

Despite the arithmetic complications involved in calculating $\mathcal{R}$, conceptually the idea reduces to an application of Taylor's theorem. On a single interval of integration $I_j$, we have that for every $x$, there exists some $\xi_x$ so that 
\begin{equation}
|F'(x)|=|F'(x_i)+(x-x_i)F''(\xi_x)|\leq |F'(x_i)|+\frac{\delta}{2} ||F''||_\infty. 
\end{equation}
We then find $\mathcal{R}$ using a combination of local and global estimates for the derivatives of the integrand corresponding to simple pointwise evaluations at the integration points and global estimates involving the neural network weights. The exact calculation of the estimates that produce the regularizer $\mathcal{R}$ are tedious and deferred to Appendix \ref{appRCalc}, with $\mathcal{R}$ given in \eqref{eqDefReg} via expressions in \eqref{eqR1} and \eqref{eqR2}. 


\section{Numerical Results} \label{sec:num_res}

In this section, we test some of the alternatives proposed in Section \ref{sec:integral_approx} to overcome the quadrature problems. Specifically, we solve the two model problems from Section \ref{sec:model_problem} using: (a) a piecewise-linear approximation of the NN, and (b) an adaptive integration method. We also solve model problem 2 using regularization methods. 

For the cases of piecewise-linear approximation and adaptive integration, we use a NN architecture composed of one hidden layer with ten neurons and a \textit{sigmoid} activation function. Moreover, we select SGD as the optimizer. In the case of regularization, we use a \textit{hyperbolic tangent} activation function with the Adam optimizer.

\subsection{Piecewise-linear approximation}

We select a piecewise-linear approximation of the NN as our approximate solution. We compute the gradients using Finite Differences and the integrals using a one-point Gaussian quadrature rule (i.e., the midpoint rule). For model problem 1, we use two different uniform partitions composed of four and ten elements, and execute $40,000$ iterations. For model problem 2, we use a uniform partition composed of ten elements and execute $200,000$ iterations.

\begin{figure}[!htp]
\centering
	\subcaptionbox{Model problem 1.\label{fig:Ritz_model_problem_1_FD_loss}}{%
	\includegraphics[width=.45\linewidth]{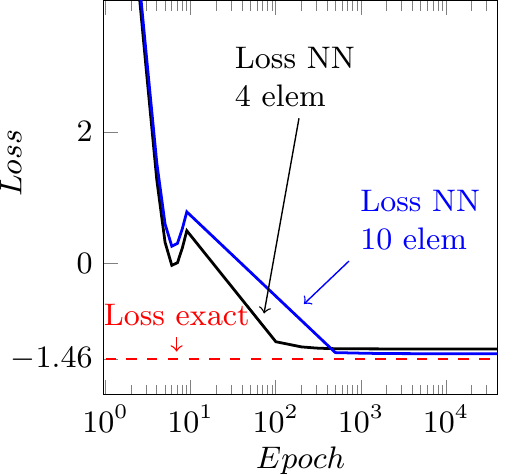}} 
	\hspace{0.5cm}
	\subcaptionbox{Model problem 2.\label{fig:Ritz_model_problem_2_FD_loss}}{%
	\includegraphics[width=.45\linewidth]{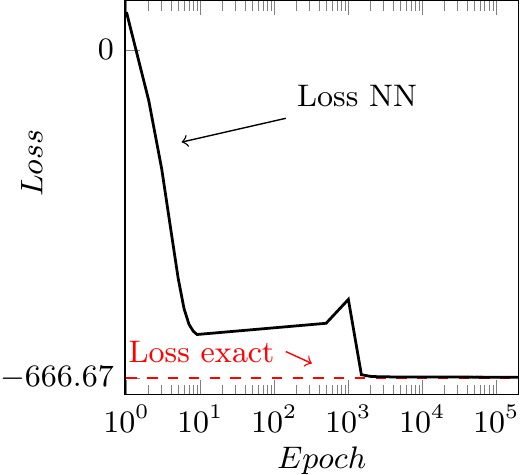}}
	\caption{Loss evolution of the training process for our two model problems when we use a piecewise-linear approximation of the NN.}
	\label{fig:Ritz_model_problem_FD_loss}
\end{figure}

\begin{figure}[!htp]
\centering
	\subcaptionbox{Model problem 1.\label{fig:Ritz_model_problem_1_linear_approx}}{%
	\includegraphics[width=.45\linewidth]{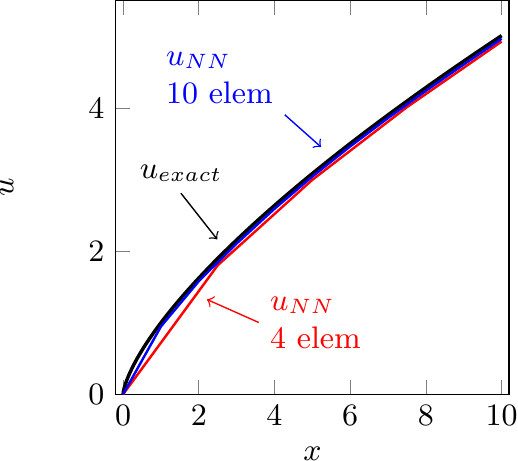}} 
	\hspace{0.5cm}
	\subcaptionbox{Model problem 2.\label{fig:Ritz_model_problem_2_linear_approx}}{%
	\includegraphics[width=.45\linewidth]{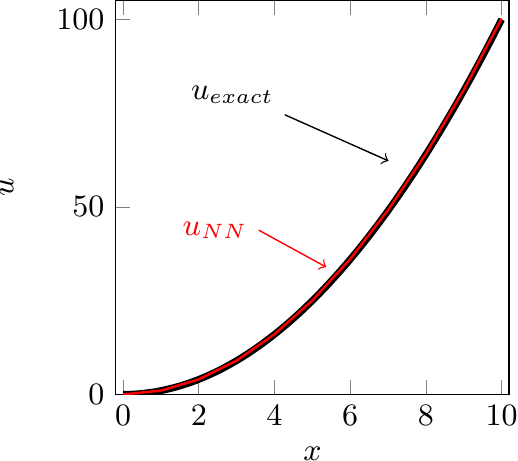}}
	\caption{Ritz method solution when we use a piecewise-linear approximation of the NN to solve the problem.}
	\label{fig:Ritz_model_problem_FD}
\end{figure}

Figures \ref{fig:Ritz_model_problem_1_FD_loss} and \ref{fig:Ritz_model_problem_2_FD_loss} show that the loss converges to the loss of the exact solution. We also observe better results as we increase the number of elements, as physically expected. Figures \ref{fig:Ritz_model_problem_1_linear_approx} and \ref{fig:Ritz_model_problem_2_linear_approx} show the corresponding solutions, which are consistent with the loss evolutions displayed in Figures \ref{fig:Ritz_model_problem_1_FD_loss} and \ref{fig:Ritz_model_problem_2_FD_loss}.

While the use of a piecewise-linear approximation overcomes the quadrature problems, the convergence is limited to $\mathcal{O}(h)$, where $h$ is the element size \cite{conv_linear_approx}. To increase this speed, it is possible to consider different piecewise-polynomial approximations, including the use of \textit{r}-adaptive algorithms \cite{r_adap}.

\subsection{Adaptive integration}

We compute the gradients using automatic differentiation and the integrals using a three-point Gaussian quadrature rule. As in previous examples, we start the training with a uniform partition of four elements for model problem 1 and of ten elements for model problem 2. We select a half-size partition of the training dataset for validation. For model problem 1, we compare the integral values of the training and validation sets (i.e., we execute Algorithm \ref{alg:adaptive}) every $100$ iterations, with an error tolerance of $0,0001$; for model problem 2, we compare every $10,000$ iterations, with an error tolerance of $10$. The tolerance is selected as a small percentage of the value of the loss.

\begin{figure}[!htp]
\centering
	\subcaptionbox{Model problem 1.\label{fig:Ritz_model_problem_1_adaptive_loss}}{%
	\includegraphics[width=.44\linewidth]{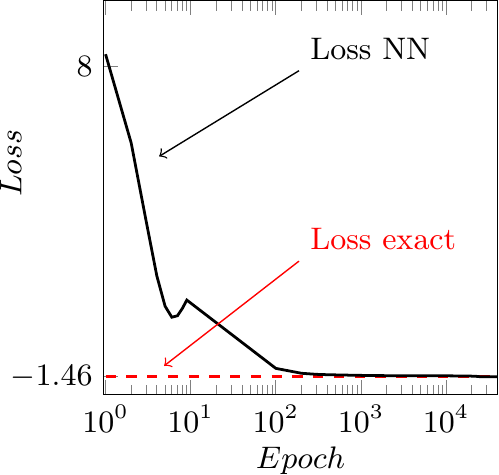}} 
	\hspace{0.5cm}
	\subcaptionbox{Model problem 2.\label{fig:Ritz_model_problem_2_adaptive_loss}}{%
	\includegraphics[width=.45\linewidth]{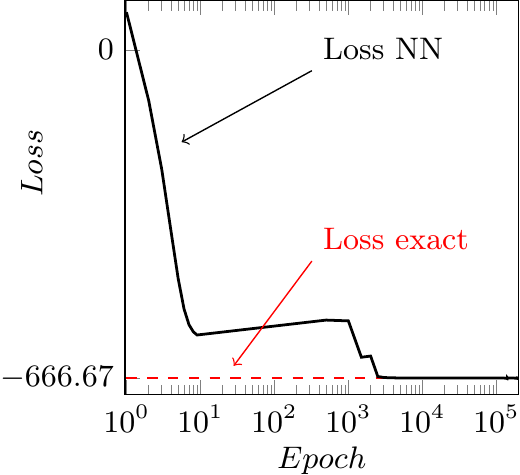}}
	\caption{Loss evolution of the training process for our two model problems when using adaptive integration.}
	\label{fig:Ritz_model_problem_adaptive_loss}
\end{figure}

\begin{figure}[!htp]
\centering
	\subcaptionbox{Model problem 1.\label{fig:Ritz_model_problem_1_adaptive}}{%
	\includegraphics[width=.45\linewidth]{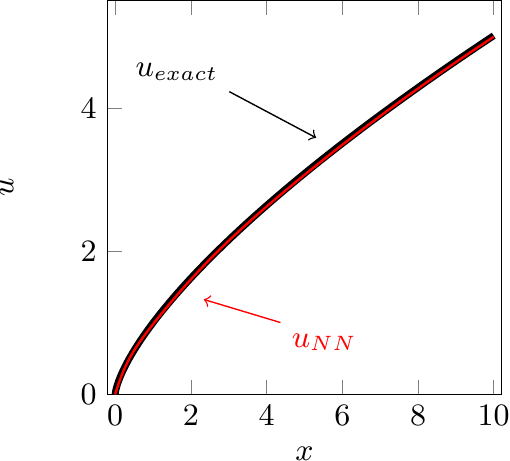}} 
	\hspace{0.5cm}
	\subcaptionbox{Model problem 2.\label{fig:Ritz_model_problem_2_adaptive}}{%
	\includegraphics[width=.45\linewidth]{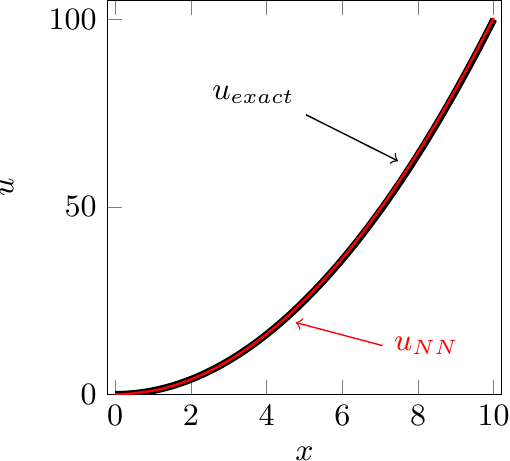}}
	\caption{Ritz method solution when using adaptive integration.}
	\label{fig:Ritz_model_problem_adaptive}
\end{figure}

Figures \ref{fig:Ritz_model_problem_1_adaptive_loss} and \ref{fig:Ritz_model_problem_2_adaptive_loss} show that the loss converges to the optimum value. As explained in Section \ref{sec:integral_approx}, the adaptive integration algorithm refines the training dataset. For model problem 1, the algorithm refines the first interval at iteration $25,000$. For model problem 2, two refinements occur in the first interval. The adaptive integration algorithm automatically selects the first interval for refinement, where overfitting was taking place in Figures \ref{fig:Ritz_model_problem_1} and \ref{fig:Ritz_model_problem_2}. Figures \ref{fig:Ritz_model_problem_1_adaptive} and \ref{fig:Ritz_model_problem_2_adaptive} show the corresponding solutions. We observe that the approximate solutions properly approximate the exact ones. 

The extrapolation of this method to higher dimensions (2D or 3D) requires higher-dimensional discretizations and quadrature rules.

\subsection{Regularization methods}

We do not apply the regularization method to model problem 1 as the method requires sufficient regularity in order to provide the necessary estimates in the calculation of $\mathcal{R}$. Since the solution is singular at $x=0$, the necessary Lipschitz bounds on the integral functional cannot be obtained within this framework. Instead, we aim to demonstrate that for problems that are sufficiently regular, our technique can avoid overfitting, and leave open the question as to how one may adapt the technique to singular problems for future work. We thus consider model problem 2. We propose the loss defined via 
\begin{equation}
\mathcal{L}(\theta)=\mathcal{\hat{F}}_{R}(u_{NN})+\mathcal{R}(\theta).
\end{equation}
Explicitly, 

\begin{equation}
\mathcal{\hat{F}}_{R}(u_{NN})=\frac{10}{N}\sum\limits_{j=1}^N\frac{1}{2}|u'_{NN}(x_j)|^2-2u_{NN}(x_j)-20u_{NN}(20),
\end{equation} 
where $x_j=\frac{10}{N}\left(i-\frac{1}{2}\right)$.

\subsubsection{Experiment 1}
We consider $N=50$ points, and a single layer network with $M=10$ neurons. We use the Adam optimizer with learning rate $10^{-2}$. We solve model problem 2 with two losses: with and without regularization. In both cases, we measure the metrics $\mathcal{L}$, $\mathcal{R}$, and $\mathcal{\hat{F}}_{R}$. For validation, we use an equidistant partition of $(0,10)$ with $49$ points, so that we still use a midpoint rule but with different integration points.

    \begin{figure*}
        \centering
        \begin{subfigure}[b]{0.475\textwidth}
            \centering
            	\includegraphics[width=.88\linewidth]{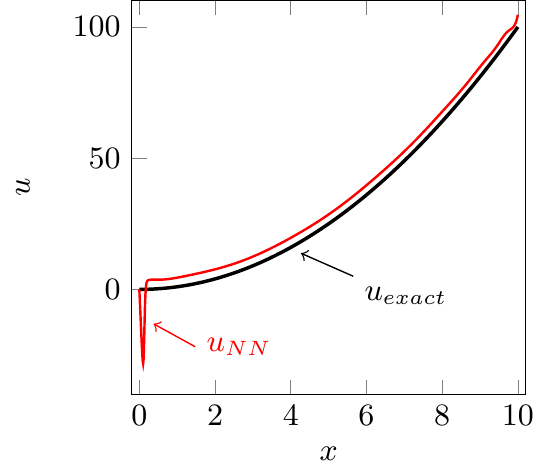}
            \caption{Exact and approximated solutions.}%
            \label{fig:1A}
        \end{subfigure}
        \hfill
        \begin{subfigure}[b]{0.475\textwidth}  
            \centering 
	\includegraphics[width=.85\linewidth]{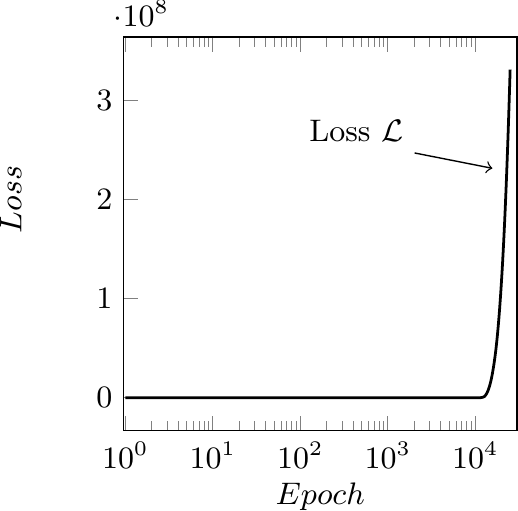}
            \caption{Evolution of $\mathcal{L}$ during training.}
            \label{fig:1B}
        \end{subfigure}
        \vskip\baselineskip
        \begin{subfigure}[b]{0.475\textwidth}   
            \centering 
	\includegraphics[width=.85\linewidth]{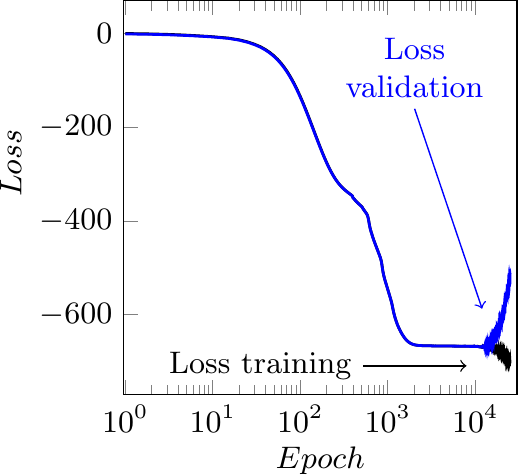}
            \caption{Evolution of $\mathcal{\hat{F}}_{R}$ during training.}
            \label{fig:1C}
        \end{subfigure}
        \hfill
        \begin{subfigure}[b]{0.475\textwidth}   
            \centering 
	\includegraphics[width=.85\linewidth]{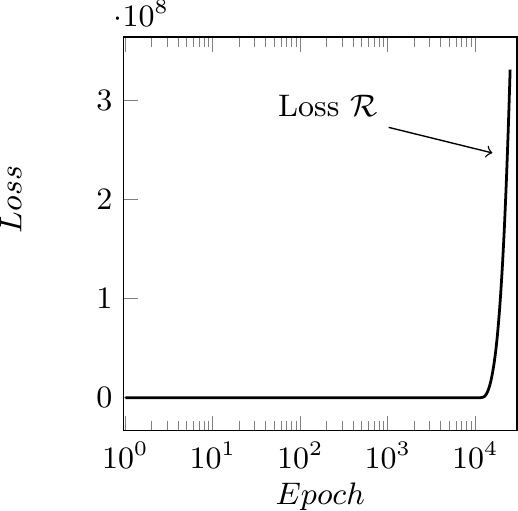}
	\caption{Evolution of $\mathcal{R}$ during training.}
            \label{fig:1D}
        \end{subfigure}
        \caption{The solution and training information for Experiment 1 without regularization.}
        \label{fig:Jamie_Imp1_NOREG}
    \end{figure*}

        \begin{figure*}
        \centering
        \begin{subfigure}[b]{0.475\textwidth}
            \centering
            	\includegraphics[width=.86\linewidth]{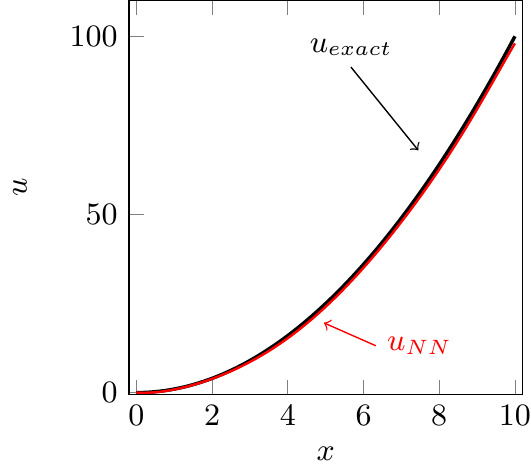}
            \caption{Exact and approximated solutions.}%
            \label{fig:2A}
        \end{subfigure}
        \hfill
        \begin{subfigure}[b]{0.475\textwidth}  
            \centering 
	\includegraphics[width=.85\linewidth]{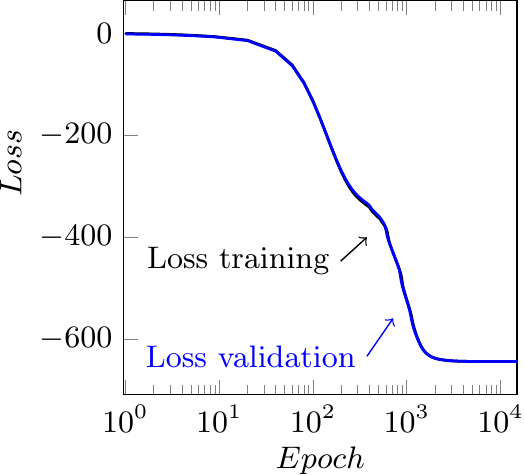}
            \caption{Evolution of $\mathcal{L}$ during training.}
            \label{fig:2B}
        \end{subfigure}
        \vskip\baselineskip
        \begin{subfigure}[b]{0.475\textwidth}   
            \centering 
	\includegraphics[width=.85\linewidth]{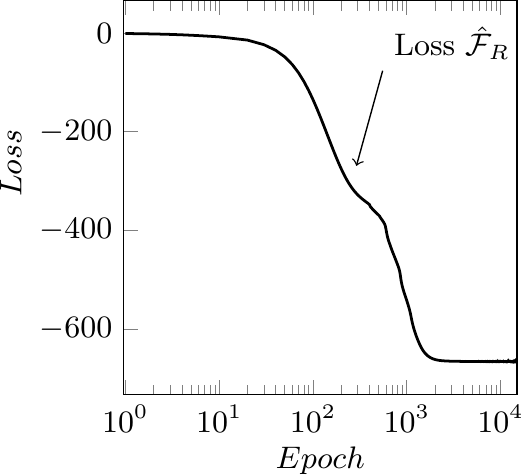}
            \caption{Evolution of $\mathcal{\hat{F}}_{R}$ during training.}
            \label{fig:2C}
        \end{subfigure}
        \hfill
        \begin{subfigure}[b]{0.475\textwidth}   
            \centering 
	\includegraphics[width=.85\linewidth]{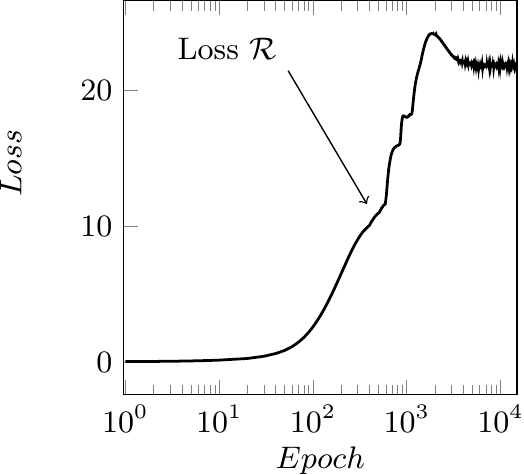}
	\caption{Evolution of $\mathcal{R}$ during training.}
            \label{fig:2D}
        \end{subfigure}
        \caption{The solution and training information for Experiment 1 with regularization.}
        \label{fig:Jamie_Imp1_REG}
    \end{figure*}

\Cref{fig:Jamie_Imp1_NOREG} shows the results without regularization. As expected, we see in \Cref{fig:1A} that the approximation is poor due to overfitting, which is most notable around $x=0$ and attained within 5000 epochs. Via the provided plots we can observe the beginning of overfitting in two distinct manners. First, we observe in \Cref{fig:1C} that the value of $\mathcal{\hat{F}}_{R}$ evaluated over the validation data begins to diverge from the value on the training data, becoming apparent at around 1000 epochs. We also see this behaviour reflected in the evolution of $\mathcal{R}$ in \Cref{fig:1D}, with its most dramatic increase beginning around the same iteration. This rapid increase also provokes an increase in $\mathcal{L}$, as seen in \Cref{fig:1B}. This also indicates that even if $\mathcal{R}$ is not used as part of the training process, its increase could be used as a metric to identify overfitting.

\Cref{fig:Jamie_Imp1_REG} describes the results with regularization and we observe a different behaviour. The approximation is generally good, and we do not see any signs of overfitting within $10^5$ epochs, as shown in \Cref{fig:2A}. In particular, the values of $\mathcal{\hat{F}}_{R}$ at the training and validation data remain consistent in \Cref{fig:2C}. Throughout \Cref{fig:Jamie_Imp1_REG} we see that within $10^5$ epochs all metrics appear to have converged to a limiting value. We obtain final values $\mathcal{L}\approx-644.22$, $\mathcal{\hat{F}}_{R}\approx -666.07$, $\mathcal{R}\approx 24.8$. We recall that the true energy of the exact solution is $\mathcal{F}_{R}(u_{exact})\approx -666.667$, which suggests the quadrature rule is accurate. Notice that in the case without regularization, before overfitting became apparent, $\mathcal{R}$ had already attained values of around 1000, which is far larger than the value of $\mathcal{R}$ at the obtained solution when regularization was used. 

\subsubsection{Experiment 2}

We now consider a smaller $N$. As we expect $\mathcal{R}$ to scale as $\frac{1}{N}$, we anticipate a more adverse effect when $N$ is small. To view this, we consider the same problem of Experiment 1, where we now select $N=20$ integration points. We consider $M=10$ neurons and minimize our problem using the Adam optimizer with a learning rate of $10^{-2}$. As before, we consider the cases with and without regularization.

    \begin{figure*}
        \centering
        \begin{subfigure}[b]{0.475\textwidth}
            \centering
            	\includegraphics[width=.87\linewidth]{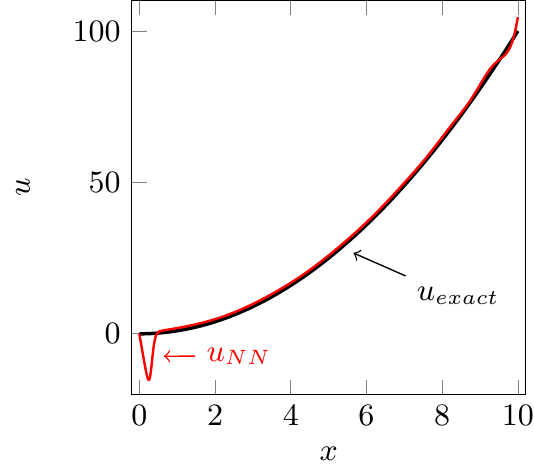}
            \caption{Exact and approximated solutions.}%
            \label{fig:3A}
        \end{subfigure}
        \hfill
        \begin{subfigure}[b]{0.475\textwidth}  
            \centering 
	\includegraphics[width=.85\linewidth]{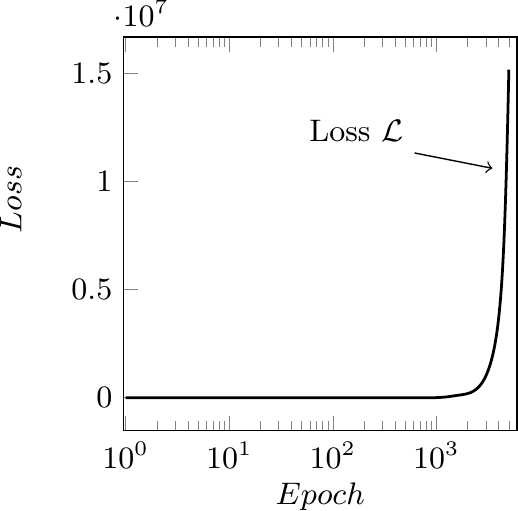}
            \caption{Evolution of $\mathcal{L}$ during training.}
            \label{fig:3B}
        \end{subfigure}
        \vskip\baselineskip
        \begin{subfigure}[b]{0.475\textwidth}   
            \centering 
	\includegraphics[width=.85\linewidth]{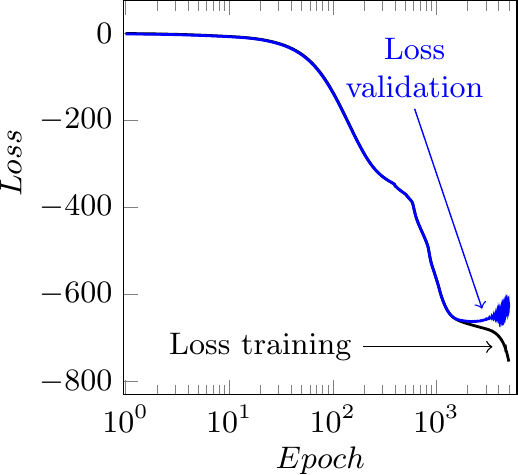}
            \caption{Evolution of $\mathcal{\hat{F}}_{R}$ during training.}
            \label{fig:3C}
        \end{subfigure}
        \hfill
        \begin{subfigure}[b]{0.475\textwidth}   
            \centering 
	\includegraphics[width=.85\linewidth]{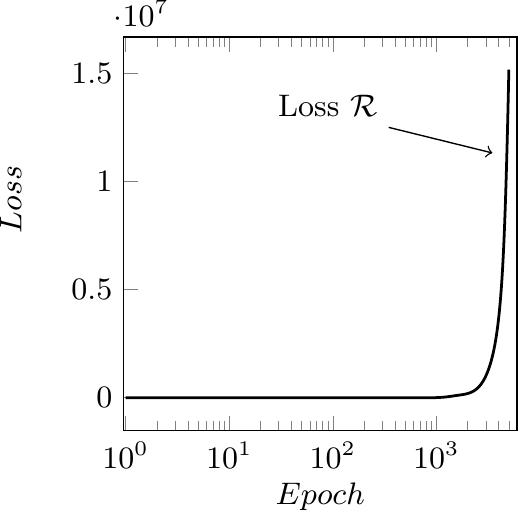}
	\caption{Evolution of $\mathcal{R}$ during training.}
            \label{fig:3D}
        \end{subfigure}
        \caption{The solution and training information for Experiment 2 without regularization.}
        \label{fig:Jamie_Imp2_NOREG}
    \end{figure*}

        \begin{figure*}
        \centering
        \begin{subfigure}[b]{0.475\textwidth}
            \centering
            	\includegraphics[width=.85\linewidth]{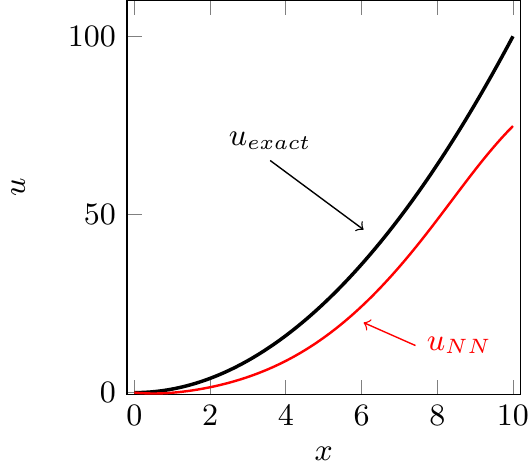}
            \caption{Exact and approximated solutions.}%
            \label{fig:4A}
        \end{subfigure}
        \hfill
        \begin{subfigure}[b]{0.475\textwidth}  
            \centering 
	\includegraphics[width=.85\linewidth]{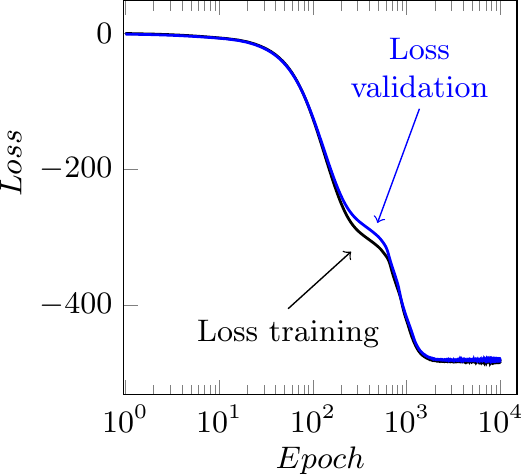}
            \caption{Evolution of $\mathcal{L}$ during training.}
            \label{fig:4B}
        \end{subfigure}
        \vskip\baselineskip
        \begin{subfigure}[b]{0.475\textwidth}   
            \centering 
	\includegraphics[width=.85\linewidth]{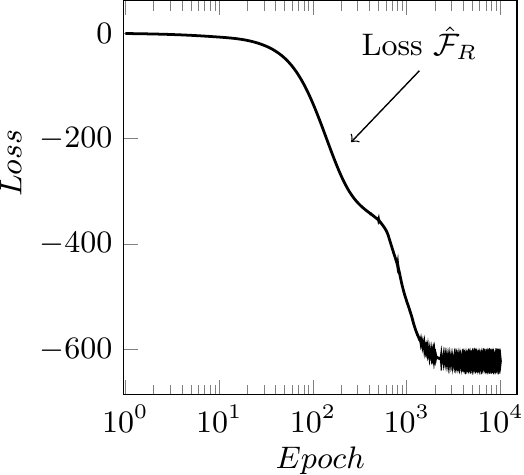}
            \caption{Evolution of $\mathcal{\hat{F}}_{R}$ during training.}
            \label{fig:4C}
        \end{subfigure}
        \hfill
        \begin{subfigure}[b]{0.475\textwidth}   
            \centering 
	\includegraphics[width=.85\linewidth]{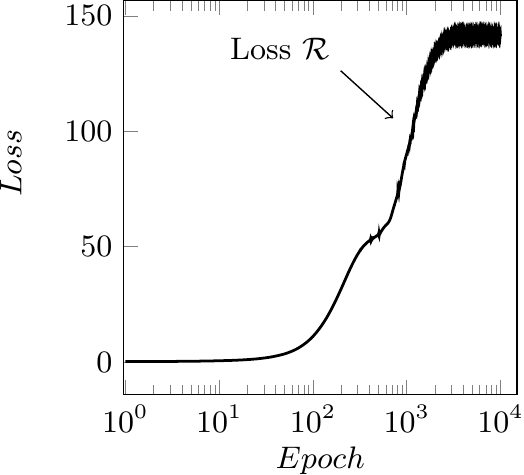}
	\caption{Evolution of $\mathcal{R}$ during training.}
            \label{fig:4D}
        \end{subfigure}
        \caption{The solution and training information for Experiment 2 with regularization.}
        \label{fig:Jamie_Imp2_REG}
    \end{figure*}
    
\Cref{fig:Jamie_Imp2_NOREG} presents the loss evolution without regularization. We observe overfitting, which is accompanied by divergence of the loss on the validation dataset, as well as a rapid increase in $\mathcal{R}$, with these features visible within 5000 epochs. 

\Cref{fig:Jamie_Imp2_REG} presents the results with regularization. We observe no signs of overfitting, with the validation and training loss remaining close in \Cref{fig:4B}. All metrics appear to have converged to a limiting value within $10^4$ epochs. However, the large value of $\mathcal{R}$ at the found solution (approximately 140) has substantially changed the optimization problem so that the obtained minimizer is far from the desired solution. The final value of $\mathcal{\hat{F}}_{R}$ is around $-622$, which is far from the desired value of $-666.67$. This experiment highlights the fact that the regularizer becomes more effective when a large number of integration points are used.


\section{Conclusions \& Future work} \label{sec:conclu}

We first illustrated how quadrature errors can destroy the quality of the approximated solution when solving PDEs using DL methods. Thus, it is crucial to select an adequate method to overcome the quadrature problems. Herein, we proposed four different alternatives: (a) Monte Carlo integration, (b) piecewise-polynomial approximations of the output of the Neural Network, (c) adaptive integration, and (d) regularization methods. We discussed the advantages and limitations of each of these methods, and we illustrated their performance via simple 1D numerical examples.

In high dimensions, Monte Carlo integration methods are the best choice. Regularizer methods are another option, but they are problem dependent and they need to be derived for each different architecture. Moreover, they require further analysis for highly nonlinear integrands. Furthermore, they are limited only to sufficiently smooth integral functionals. In addition, more complex NN architectures (which should be needed in higher dimension) will hinder the derivation of $\mathcal{R}$.

In low dimensions (three or below), Monte Carlo integration is not competitive because of its low convergence speed. In these cases, adaptive integration exhibits faster convergence. In the cases of piecewise-linear approximation and regularizers, we are also able to overcome the quadrature problems, but the convergence speed is often slower than with adaptive integration.

Possible future research lines of this work are: (a) to implement adaptive integration method in 2D and 3D, (b) to apply \textit{r-adaptivity} methods with piecewise-polynomial approximations, and (c) to implement regularizers for high dimensional problems involving more complex NN architectures.

\section{Acknowledgments}

This work has received funding from: the European Union's Horizon 2020 research and innovation program under the Marie Sklodowska-Curie grant agreement No 777778 (MATHROCKS); the European Regional Development Fund (ERDF) through the Interreg V-A Spain-France-Andorra program POCTEFA 2014-2020 Project PIXIL (EFA362/19); the Spanish Ministry of Science and Innovation projects with references PID2019-108111RB-I00 (FEDER/AEI), PDC 2021-121093-I00, and PID2020-114189RB-I00 and the ``BCAM Severo Ochoa" accrediation of excellence (SEV-2017-0718); and the Basque Government through the three Elkartek projects 3KIA (KK-2020/00049), EXPERTIA (KK-2021/000 48), and SIGZE (KK-2021/00095), the Consolidated Research Group MATHMODE (IT1294-19) given by the Department of Education, and the BERC 2018-2021 program.


\appendix

\section{Estimation of the regularizer}
\label{appRCalc}

Following the heuristic of \Cref{subsecHeuristic}, we derive an expression for $\mathcal{R}$ that may be used as a regularizer. The necessary steps are:
\begin{enumerate}
\item Using the chain rule, we find global upper bounds for the derivatives of a simple neural network in terms of the weights. 
\item Via Taylor's theorem with remainder and using the global estimates for simple networks, we find local estimates of the derivatives of NNs with a cutoff function to ensure a homogeneous Dirichlet condition. 
\item Using local estimates for the derivatives of a NN, we find local Lipschitz estimates for integrands corresponding to the Ritz method.
\end{enumerate}

We tackle each of these estimations in the following subsections. 

\subsection{Global estimates for derivatives of a single layer network}
Let $\hat{u}_{NN}$ be a single layer neural network. We write it in the form 
\begin{equation}\label{eqDefNet}
\hat{u}_{NN}(x)=b^1+\sum\limits_{i=1}^MA_i^1\sigma(A_i^0x+b^0_i).
\end{equation}
for weights $\theta=(b^1,A^1,b^0,A^0)$ and an activation function $\sigma$. We assume that $\sigma$ has globally bounded derivatives, so that $||\sigma^{(n)}||_\infty$ is finite for every $n=0,1,2...$. 

A simple application of the triangle quality and that $\sigma$ is bounded gives that 
\begin{equation}
\begin{split}
|\hat{u}_{NN}(x)|\leq &|b^1|+\sum\limits_{i=1}^M |A_i^1\sigma(A_i^0x+b^0_i)|\\
\leq & |b^1|+\sum\limits_{i=1}^M|A_i^1|||\sigma||_\infty\\
\end{split}
\end{equation}

The derivatives of $\hat{u}_{NN}$ are 
\begin{equation}
\hat{u}^{(n)}_{NN}(x)=\sum\limits_{i=1}^MA_i^1(A_i^0)^n\sigma^{(n)}(A_i^0x+b^0_i),
\end{equation}
for $n\geq 1$. Similarly, this gives the immediate estimate that 
\begin{equation}
\left|\hat{u}^{(n)}_{NN}(x)\right|\leq \sum\limits_{i=1}^M|A_i^1||A_i^0|^n||\sigma^{(n)}||_\infty.
\end{equation}

Thus, we define the global upper bounding function $\mathcal{R}^1(\theta;n)$ by 

\begin{equation}\label{eqR1}
\mathcal{R}^1(\theta;n)=\left\{\begin{array}{c c}
|b^1|+||\sigma||_\infty\sum\limits_{i=1}^M|A_i^1| & n=0,\\
||\sigma^{(n)}||_\infty\sum\limits_{i=1}^M|A_i^1||A_i^0|^n & n\geq 1,
\end{array}\right.
\end{equation}
which gives the global upper bound 
\begin{equation}
\left|\hat{u}_{NN}^{(n)}(x)\right|\leq \mathcal{R}^1(\theta;n)
\end{equation}
for every $x$.

\subsection{Local derivative estimation of a single layer neural network with cutoff function}

Next, we consider the commonly considered case where our neural network admits a cutoff function to ensure a Dirichlet boundary condition, and turn to the case of {\it local} estimation. Explicitly, we take $u_{NN}(x)=\hat{u}_{NN}(x)\phi(x)$, where $\phi$ is zero at the homogeneous Dirichlet condition, and $\hat{u}_{NN}$ is as in \eqref{eqDefNet}. We presume that $\phi$ has bounded derivatives; i.e. $||\phi^{(k)}||_\infty$ is finite for $k=0,1,...$. Let $x$ be in the interval $I_j=\left(x_j-\frac{\delta}{2},x_j+\frac{\delta}{2}\right)$. Then, we have

\begin{equation}
\begin{split}
u^{(n)}_{NN}(x)=&u^{(n)}_{NN}(x_j)+(x-x_j)u^{(n+1)}_{NN}(\xi_x)\\
=&u^{(n)}_{NN}(x_j)+(x-x_j)\sum\limits_{k=0}^{n+1}\binom{n+1}{k}u^{(k)}_{NN}(\xi_x)\phi^{(n+1-k)}(\xi_x)
\end{split}
\end{equation}
via Taylor's theorem with remainder and the product rule for higher order derivatives. Employing the global upper bounds $\mathcal{R}^1(\theta;n)$ from \eqref{eqR1}, we have that for $x\in I_j$, 
\begin{equation}
\begin{split}
\left|u^{(n)}_{NN}(x)\right|\leq & \left|u^{(n)}_{NN}(x_j)\right|+\frac{\delta}{2}\sum\limits_{k=0}^{n+1}\binom{n+1}{k}|u^{(k)}_{NN}(\xi_x)||\phi^{(n+1-k)}(\xi_x)|\\
\leq &  \left|u^{(n)}_{NN}(x_j)\right|+\frac{\delta}{2}\sum\limits_{k=0}^{n+1}\binom{n+1}{k}\mathcal{R}^1(\theta;k)||\phi^{(n+1-k)}||_\infty.
\end{split}
\end{equation}

Thus, we define the second intermediate regularizer as 

\begin{equation}\label{eqR2}
\mathcal{R}^2(\theta;I_j,n)=\left|u^{(n)}_{NN}(x_j)\right|+\frac{\delta}{2}\sum\limits_{k=0}^{n+1}\binom{n+1}{k}\mathcal{R}^1(\theta;k)||\phi^{(n+1-k)}||_\infty,
\end{equation}
giving the estimate that for all $x\in I_j$,
\begin{equation}
\left|u^{(n)}_{NN}(x)\right|\leq \mathcal{R}^2(\theta;I_j,n).
\end{equation}
We note that via automatic differentiation, $u^{(n)}$ may be evaluated at the training data $x_j$.

\subsection{Application to integral errors}

Our aim is to estimate the error in the integral functional
\begin{equation}
\mathcal{F}(u_{NN})=\int_a^b \frac{1}{2}|u_{NN}'(x)|^2-f(x)u_{NN}(x)\,dx-g(a)u_{NN}(a)-g(b)u_{NN}(b),
\end{equation}
when approximated by a simple quadrature rule 
\begin{equation}\label{eqQuadRuleAp}
\sum\limits_{j=1}^N \left(\frac{1}{2}|u_{NN}'(x_j)|^2-f(x_j)u_{NN}(x_j)\right)\delta -g(b)u_{NN}(b)-g(a)u_{NN}(a).
\end{equation}
The boundary terms can be calculated in one dimension without quadrature error and thus we ignore their contribution. We estimate the error for the quadrature rule \label{eqQuadRuleAp} by obtaining Lipschitz bounds of the integrand via

\begin{equation}
\begin{split}
&\left|\frac{d}{dx}\left(\frac{1}{2}|u_{NN}'(x)|^2-f(x)u_{NN}(x)\right)\right|\\
=& \left|u'_{NN}(x)u''_{NN}(x)-f'(x)u_{NN}(x)-f(x)u'_{NN}(x)\right|\\
\leq & |u'_{NN}(x)||u''_{NN}(x)|+|f'(x)||u_{NN}(x)|+|f(x)||u'_{NN}(x)|.
\end{split}
\end{equation}

Estimating the (local) Lipschitz constant of the integrand reduces to estimating (locally) various derivatives of $u_{NN}$. For $x\in I_j$, we estimate the Lipschitz constant of the integrand via 
\begin{equation}
\left|\frac{d}{dx}\left(\frac{1}{2}|u_{NN}'(x)|^2-f(x)u_{NN}(x)\right)\right|\leq \mathcal{R}^3(\theta;I_j),
\end{equation}
where the regularizer $\mathcal{R}^3(\theta,I_j)$ is given by 
\begin{equation}
\mathcal{R}^3(\theta;I_j)=\Big(\mathcal{R}^2(\theta;I_j,1)\mathcal{R}^2(\theta;I_j,2)+||f||_\infty \mathcal{R}^2(\theta;I_j,1)+||f'||_\infty \mathcal{R}^2(\theta;I_j,0)\Big).
\end{equation}

We define the final regularizer $\mathcal{R}$ by 
\begin{equation}\label{eqDefReg}
\mathcal{R}(\theta)=\frac{\delta^2}{4}\sum\limits_{j=1}^N \mathcal{R}^3(\theta;I_j),
\end{equation}
which following \eqref{eqMidpError} gives the estimate 
\begin{equation}\label{eqError}\begin{split}
&|\mathcal{F}(u_{NN})-\hat{\mathcal{F}}(u_{NN})|\\
=&\left|\int_a^b \frac{1}{2}|u_{NN}'(x)|^2-f(x)u_{NN}(x)\,dx - \sum\limits_{j=1}^N \left(\frac{1}{2}|u_{NN}'(x_j)|^2-f(x_j)u_{NN}(x_j)\right)\delta\right|\\
\leq &\mathcal{R}(\theta). 
\end{split}
\end{equation}

\bibliography{biblio}

\end{document}